# Generalized Kantorovich-type theorem for the Fixed Slope Iterations

Andrei Dubin[1]
ITEP, 117218, B.Cheremushkinskaya 25, Moscow, Russia

## Abstract

The extended modification of the Newton method is considered when the inverse of the derivative $F'(x)$ (of the operator $F(x)$ in the equation $F(x) = 0$) is replaced by an invertible bounded $x$-independent operator $B$. The continuity assumption is relaxed to the requirement that $F(x)$ is continuously Frechet-differentiable. The Kantorovich majorization technique is adapted to formulate and prove the corresponding generalization of the Kantorovich theorem originally stated for the standard modified Newton method (with $B^{-1} = F'(x_0)$) when $F'(x)$ is Lipschitz continuous. If $B^{-1} = F'(x_0)$, the generalized theorem is shown to extend the existing one due to Argyros. For a generic $B$ and a Holder continuous $F'(x)$, the proposed theorem leads to a weaker condition of the semilocal convergence, larger uniqueness domain and finer error bounds compared to the previous results of Ahues and Argyros.

## 1   Introduction

Consider the Newton-like method defined by the recurrence relation

$$x_{k+1} = x_k - BF(x_k) = T(x_k) \quad , \quad k \geq 0, \qquad (1\text{-}1)$$

where $F(x): D_F \subseteq X \to Y$ is a continuously Frechet-differentiable operator which defines a mapping from an open subset $D_F$ of a Banach space $X$ into a given Banach space $Y$. In turn, an invertible bounded linear operator $B$ may be viewed as an approximation to $(F'(x_0))^{-1}$. The $B = A^{-1}$ implementation of ( 1-1 ) was considered under different names by many authors as such simpler alternative to the Newton method that does not require the numerical inversion of $F'(x_k)$ at each point $x_k$. Recently such implementation attracted renewed attention under the name of the fixed slope iterations (FSI) method applied to the case when the derivative operator $F'(x)$ is Lipschitz [2] and Holder [3], [4] continuous.

   The purpose of the present paper is to formulate and prove a theorem which generalizes the Kantorovich theorem [1] for the modified Newton method (MNM) with $B^{-1} = F'(x_0)$. The proposed extension deals with the general case of the continuously differentiable (rather than Lipschitz continuous as in [1]) derivative $F'(x)$ that relaxes the continuity requirements of [ 2 ]-[ 4 ]. One notes also that *no* inversion of either $B$ or $F'(x_0)$ is required by the computations. In consequence, the considered numerical analysis of the method ( 1-1 ) is not restrained when $F'(x_0)$ or $B$ is ill-conditioned.

---

[1] E-mail: dubin@itep.ru

When reduced to the MNM case $B^{-1} = F'(x_0)$, the proposed implementation of the FSI method may be viewed as such extension of the approach [ 5 ] that introduces a more workable form of the convergence condition and a larger uniqueness ball. In addition, it is proved that thus generated Newton-like sequence $x_k$ converges to one and the same limiting point $x_*$ for *any* initial guess $\tilde{x}_0$ in the extended uniqueness domain. It is also noted that the uniqueness domain of the FSI method is generically wider than the one of the contraction mapping method.

In the specific case of the Holder continuous $F'(x)$ and generic invertible $B \in L(Y, X)$, the new theorem improves the FSI conditions of the semilocal convergence and uniqueness as well as the error bounds introduced in [ 2 ]-[ 4 ]. The latter FSI technique introduces the parameter $\delta \in [0,1[$ to measure how close the operators $A = B^{-1}$ and $F'(x_0)$ are. In the limit $\delta \to +0$, the convergence condition [ 3 ], [ 4 ] is more restrictive than the existing MNM convergence condition (e.g., see [ 6 ]) obtained when $F'(x)$ is Holder continuous. Moreover, when $F'(x)$ is Lipschitz continuous, the limit $\delta \to +0$ of the condition [ 2 ]-[ 4 ] is suboptimal compared to the well-known MNM Kantorovich condition [ 1 ]. On the contrary, the proposed new theorem in the limit $\delta \to +0$ does reproduce the latter two existing MNM results.

## 2 Generalized Kantorovich-type theorem

### 2.1 Adaptation of the Kantorovich majorization technique

To formulate and prove the convergence theorem, one is to adapt the following two theorems proved in Section 2 of Chapter XVIII of [ 1 ].

**Theorem 1**. Let $T(x): D_F \subseteq X \to Y$ be continuously Frechet-differentiable operator and suppose that there exists such continuously differentiable scalar function $\phi(\upsilon):[0,r] \to \mathbb{R}$ that

$$\|T'(x)\| \le \phi'(\upsilon) \quad , \quad \|x - x_0\| \le \upsilon \le r \quad , \quad \bar{B}(x_0, r) \subset D_F, \qquad (2\text{-}1)$$

$$\|Tx_0 - x_0\| \le \phi(0) = \eta, \qquad (2\text{-}2)$$

where $\bar{B}(x_0, r)$ denotes the closed ball of the center $x_0$ and the radius $r$. Presume also that the fixed point equation

$$\upsilon = \phi(\upsilon) \qquad (2\text{-}3)$$

has a minimal solution $\upsilon_* \in [0, r]$. Then, the sequence $\{x_k = T^k x_0\}: x_{k+1} = T(x_k)$ is well-defined, remains in $\bar{B}(x_0, r)$ and converges to a solution $x_* \in \bar{B}(x_0, r)$ of the equation fixed point equation $x = T(x)$. Moreover, the majorization conditions

$$\|x_{k+1} - x_k\| \le \upsilon_{k+1} - \upsilon_k \quad , \quad \|x_* - x_k\| \le \upsilon_* - \upsilon_k \quad , \quad \forall k \ge 0, \qquad (2\text{-}4)$$

are valid. The scalar majorizing sequence $\{\upsilon_k\}$, being generated by the recurrence relation $\upsilon_{k+1} = \phi(\upsilon_k)$ for $\forall k \geq 0$ with $\upsilon_0 = 0$, is non-decreasing, remains in $[0, \upsilon_*]$ and converges to a minimal solution $\upsilon_*$ of ( 2-3 ) in $[0, r]$.

**Theorem 2**. Let the conditions of the previous theorem be fulfilled and suppose that $\phi(r) \leq r$ while ( 2-3 ) has a unique solution in $[0, r]$. Then, the equation $x = T(x)$ has a unique root $x_*$ in $\bar{B}(x_0, r)$ and the sequence $\{\tilde{x}_k \equiv T^k \tilde{x}_0\}$ converges to this root starting from any $\tilde{x}_0 \in \bar{B}(x_0, r)$.

In the FSI case, the continuous non-decreasing function $\phi(\upsilon) \geq 0$ assumes the form

$$\phi(\upsilon) = \eta + \nu \cdot \upsilon + \int_0^\upsilon \bar{\omega}_0(l) dl \quad , \quad \bar{\omega}_0(\upsilon) = \omega_B(\upsilon) - \nu \quad , \quad \omega_B(0) - \nu = \bar{\omega}_0(0) = 0. \quad (\textbf{2-5})$$

Here, $\phi'(r)$ is to be equated with the continuous *non-decreasing* function $\omega_B(r)$ which, in compliance with ( 2-1 ), is to implement the affine-invariant upper bound

$$\|T'(x)\| = \|BF'(x) - 1\| \leq \omega_B(\|x - x_0\|) \quad , \quad \omega_B(\upsilon) = \phi'(\upsilon) \quad , \quad \|T'(x_0)\| \leq \nu, \quad (\textbf{2-6})$$

where $\|x - x_0\| \leq R$, $\bar{B}(x_0, R) \subset D_F$. The monotonicity of $\omega_B(\upsilon)$ implies that $\|T'(x)\| \leq \omega_B(\upsilon)$ for $\|x - x_0\| \leq \upsilon$ which, in turn, justifies the identification $\omega_B(\upsilon) = \phi'(\upsilon)$. Note that the continuity measure $\omega_B(r)$ is *centered* at $x_0$, i.e. depends on the relative distance $\|x - x_0\|$ of $x$ from the *stationary* point identified with $x_0$. The MNM case is reproduced when $\nu \to 0$.

In general, thus defined quantity $\omega_B(\upsilon)$ may be identified with a majorant of the global Lipschitz constant of the (continuously Frechet-differentiable) operator $T(x)$ on $\bar{B}(x_0, \upsilon)$ while the optimal choice of $\omega_B(\upsilon)$ is the supremum of $\|T'(x)\|$ evaluated over $\bar{B}(x_0, \upsilon)$:

$$\|T(x) - T(y)\| \leq \omega_B(\upsilon) \|x - y\| \quad , \quad \sup_{x \in \bar{B}(x_0, \upsilon)} \|T'(x)\| \leq \omega_B(\upsilon), \quad (\textbf{2-7})$$

where $x, y \in \bar{B}(x_0, \upsilon)$ and the matrix norm is presumed to be consistent with the vector norm.

## 2.2 Statement of the FSI convergence theorem

**Condition set** $\Upsilon$. Let $F(x): D_F \subseteq X \to Y$ and $B \in L(Y, X)$ be a continuously Frechet-differentiable operator and an invertible bounded linear operator respectively. Assume that the condition ( 2-6 ) (with $x \in \bar{B}(x_0, R) \subset D_F$) and the restriction $\eta \geq \|BF(x_0)\| = \|x_1 - x_0\|$ are valid for some constants $\eta > 0$, $0 \leq \nu < 1$, $R > 0$ and a continuous non-decreasing function $\omega_B(\upsilon) \geq 0$. Let the function ( 2-5 ) specify both ( 2-3 ) and the relation $\upsilon_{k+1} = \phi(\upsilon_k)$ which, being supplemented by the initial condition $\upsilon_0 = 0$, determines the scalar sequence $\{\upsilon_k\}$. Given $\nu = \omega_B(0) < 1$, let the constraint A be valid: $\phi(\gamma_*) \leq \gamma_*$ where $\gamma_* = \sup_{\gamma \in ]0, R]} \{\gamma \mid \omega_B(\gamma) < 1\}$ so that $\gamma_* \equiv \gamma_*(R) \in ]0, R]$. Finally, introduce the quantity $\lambda_*(R) \equiv \lambda_*$ in the following way. Consider the set $\Lambda$ which includes all those $\lambda \in ]\upsilon_*, R]$ that $\phi(\upsilon) - \upsilon < 0$ for $\forall \upsilon \in ]\upsilon_*, \lambda[$. Let

$\lambda_* = \sup_{\lambda \in \Lambda}[\lambda] > \upsilon_*$ ($\lambda_*$ is the supremum of $\lambda$ in $\Lambda$) if $\Lambda \neq \emptyset$ and $\lambda_* = \upsilon_*$ when $\Lambda = \emptyset$ so that $\lambda_* \in [\upsilon_*, R]$. In turn, define the conditions $B_1 = [\lambda_* = \upsilon_*] \cup [\{\lambda_* > \upsilon_*\} \cap \{\phi(\lambda_*) < \lambda_*\}]$ and $B_2 = [\{\lambda_* > \upsilon_*\} \cap \{\phi(\lambda_*) = \lambda_*\}]$.

**FSI Theorem.** Under the set $\Upsilon$ of the conditions, the sequence $\{x_k\} \equiv \{T^k x_0\}$ is well-defined, remains in $\bar{B}(x_0, \upsilon_*)$ and converges to a solution $x_*$ of the equation $BF(x) = x - T(x) = 0$. The estimates ( 2-4 ) are valid where the scalar sequence $\{\upsilon_k\}$, being non-decreasing, remains in $[0, \upsilon_*]$ and converges to a minimal solution $\upsilon_*$ of ( 2-3 ) in $[0, R]$ with $\upsilon_* \leq \gamma_*$. The solution $x_*$ is unique respectively in the closed ball $\bar{B}(x_0, \lambda_*)$ and in the open $B(x_0, \lambda_*)$ ball under $B_1$ and $B_2$. The sequence $\{\tilde{x}_k \equiv T^k \tilde{x}_0\}$ converges to $x_*$ for any $\tilde{x}_0$ in the uniqueness ball.

## 2.3 Proof of the FSI theorem

First, consider the convergence part which adapts the statements of the Theorem 1 to the FSI case. In view of ( 2-5 ) and ( 2-6 ), the conditions of the Theorem 1 are optimally fulfilled provided that the constraint A is the necessary and sufficient condition of the existence of a minimal solution $\upsilon_* \in [0, R]$ of ( 2-3 ) implemented according to ( 2-5 ). To verify it, observe first that, if $g(R) = \phi(R) - R > 0$, the *only* possible way for ( 2-3 ) to have a solution in $[0, R]$ is to maintain that the minimal value of the function $g(\upsilon)$ on $[0, R]$ is nonpositive. As $g(0) = \eta > 0$ and $g(R) > 0$, it implies that the equation $g'(\upsilon) = \omega_B(\upsilon) - 1 = 0$ has such minimal solution $\bar{r}_* \in [0, R[$ that $\phi(\bar{r}_*) \leq \bar{r}_*$ which leads to the identification $\gamma_* = \bar{r}_* \geq \upsilon_*$. In the remaining case of $g(R) \leq 0$ (with $g(0) > 0$), a minimal solution $\upsilon_* \in [0, R]$ exists by continuity of $g(\upsilon)$ and there are the two options: the equation $g'(\upsilon) = 0$ either has a minimal solution $\bar{r}_* \in [0, R]$ or it doesn't. As previously, the first option leads to $\phi(\bar{r}_*) \leq \bar{r}_*$ and $\gamma_* = \bar{r}_* \geq \upsilon_*$. As for the second option, the relation $g'(0) = -(1-\nu)$ implies that $g'(\upsilon) < 0$ for $\forall \upsilon \in [0, R]$, i.e., $\gamma_* = R$ so that $\phi(\gamma_*) \leq \gamma_*$.

To verify the uniqueness part, let us show that the conditions of the Theorem 2 are optimally satisfied. To this end, observe that $\phi(r) \leq r$ for $\forall r \in [\upsilon_*, \lambda_*]$ while $]\upsilon_*, \lambda_*[$ is the *largest* connected interval (if nonempty) of the form $]\upsilon_*, \lambda[$ which is free of solutions of ( 2-3 ). Under $B_1$, the minimal solution $\upsilon_*$ of ( 2-3 ) is unique in $[0, \lambda_*]$ that is obvious if $\upsilon_* = \lambda_*$. If $\lambda_* > \upsilon_*$, one is to use additionally that the condition $\phi(\lambda_*) < \lambda_*$ guarantees that $\lambda_*$ is not a solution of ( 2-3 ). Under $B_2$, one similarly obtains that the minimal solution $\upsilon_*$ is unique in $[0, r]$ for $\forall r \in [\upsilon_*, \lambda_*[$ because $\lambda_*$ is a solution of ( 2-3 ) in this case. Since the option $\phi(\lambda_*) > \lambda_*$ is forbidden by construction of $\lambda_*$, it completes the verification.

## 3 Application to the case of the Holder continuous operator $F'(x)$

When $F'(x)$ is Holder continuous, the upper bound ( 2-6 ) is to be introduced in the form

$\omega_B(\upsilon) = \bar{\omega}_0(\upsilon) + \nu \leq \omega_0(\upsilon) + \nu$ with $\omega_0(r) = l_0 r^\alpha$ where $\alpha \in ]0,1]$, $l_0$ is the center-Holder constant [4] and the non-decreasing continuous function $\omega_0(\upsilon)$ is in general defined as the upper bound

$$\| B(F'(x) - F'(x_0)) \| \leq \omega_0(\| x - x_0 \|) \tag{3-1}$$

to be compared with (2-6). In this case, the implementation (2-5) of (2-3) is reduced to the equation $g(\upsilon) = (1+\alpha)^{-1} l_0 \upsilon^{\alpha+1} - (1-\nu) \cdot \upsilon + \eta = 0$. The unique (global) minimum of $g(\upsilon)$ is to be equated with the solution $\bar{r}_* = [(1-\nu)/l_0]^{1/\alpha}$ of the equation $g'(\upsilon) = l_0 \upsilon^\alpha - (1-\nu) = 0$. As $g(0) = \eta > 0$, the necessary and sufficient condition for the equation $g(\upsilon) = 0$ to have a solution $\upsilon_* > 0$ is that $g(\bar{r}_*) \leq 0$ in which case the minimal solution complies with the restriction $\upsilon_* \in ]0, \bar{r}_*]$. In sum, it leads to the condition

$$l_0 \eta^\alpha \leq (1-\nu)^{\alpha+1} [\alpha/(1+\alpha)]^\alpha = l_0 \eta_{\max}^\alpha \tag{3-2}$$

that in the $\nu = 0$ case is reduced to the MNM semilocal convergence condition [6] (the equation (2.4.27)) under the Holder continuity of $F'(x)$. Given (3-2), the solution $x_*$ exists if $R \geq \upsilon_*$. It is unique in $B(x_0, \upsilon_{**}) \cap \bar{B}(x_0, R)$ and $\bar{B}(x_0, \upsilon_{**}) \cap \bar{B}(x_0, R)$ respectively when $\eta < \eta_{\max}$ and $\eta = \eta_{\max}$ where $\upsilon_{**} \geq \upsilon_*$ denotes the maximal solution of $g(\upsilon) = 0$ with $\upsilon_{**} = \upsilon_*$ if $\eta = \eta_{\max}$.

When $F'(x)$ is Lipschitz continuous ($\alpha = 1$) and $\nu = 0$, (3-2) reproduces the improved variant $2l_0 \eta \leq 1$ of the famous Kantorovich condition [1] where the Lipschitz constant $l$ is replaced [6] by the center-Lipschitz constant $l_0 \leq l$. The corresponding option of the above uniqueness condition is identical to the one in the Kantorovich theorem [1].

## 4 Comparison with the existing results

When $\nu = 0$ and $B^{-1} = F'(x_0)$, the proposed convergence condition improves the one formulated in [5] for the MNM case. The condition of [5] takes the form of the condition introduced by the Theorem 2 stated in subsection 2.1. The latter form is not fully operational because it is not directly applicable when either $\phi(r) > r$ or the equation (2-3) has multiple solutions in $[0, r]$. It is the proposed advanced formulation in terms of $\gamma_*(R)$ and $\lambda_*(R)$ that overcomes this shortage. Also, the derivation [5] of the uniqueness ball implicitly restricts it by $B(x_0, \bar{r}_*)$ where $\bar{r}_* < \lambda_*$ (if $\eta < \eta_{\max}$) denotes the minimal solution of $g'(\upsilon) = \omega_0(\upsilon) - 1 = 0$ while the non-decreasing function $\omega_0(r)$ is defined by (3-1) with $B^{-1} = F'(x_0)$. Due to (2-7), the latter restriction of the convergence ball is linked to the applicability of the contraction mapping principle.

Next, let us compare the results of the previous section with the corresponding statements of [4] where the improvement of [3] is developed via the centered reformulation in terms of the center-Holder (rather than Holder) constant $l_0 \leq l$. The considered in the previous section implementation $g(\upsilon) = \phi(\upsilon) - \upsilon = 0$ of the equation (2-3) is replaced in [4] by the suboptimal equation $f(\upsilon) = l_0 \upsilon^{\alpha+1} - (1-\delta) \cdot \upsilon + \eta = 0$ where $\delta = \nu$ (when $A^{-1} = B$) and $f(\upsilon) > g(\upsilon)$ for $\alpha > 0$. As a result the condition (3-2) is replaced by its more restrictive counterpart

$l_0\eta^\alpha \leq (1-\nu)^{\alpha+1}[\alpha/(1+\alpha)]^\alpha(1+\alpha)^{-1}$. In the $\nu = \alpha = 1$ case of a Lipschitz continuous $F'(x)$, the latter counterpart reduces to the condition $4l_0\eta \leq 1$ ([2]) that leads to the twice smaller value of $\eta_{max}$ than the one implied by centered Kantorovich condition $2l_0\eta \leq 1$. In turn, the uniqueness is verified only in the smaller ball $\bar{B}(x_0, r_*)$ where $r_*$ is the minimal solution of $f(\upsilon) = 0$ with $r_* < \upsilon_* \leq \upsilon_{**}$ and $\upsilon_*, \upsilon_{**}$ denote the minimal and maximal solution of the equation $g(\upsilon) = 0$ considered in the previous section. Note also that, given the invertibility of $B$ and the restriction $\nu = \delta < 1$, the operator inversion lemma [6] implies that $F'(x_0)$ is invertible as well.

Finally, the convergence condition is formulated in [2]-[5] only for the Newton-like sequence starting from the stationary point $x_0$ rather than from a generic point $\tilde{x}_0$ in the appropriate convergence ball.

# 5 References

[1] L.V. Kantorovich, G.P. Akilov, Functional Analysis, Pergamon Press, Oxford, 1982.

[2] M. Ahues, A note on perturbed fixed slope iterations, Applied Mathematics Letters, 18, N 4 (2005), 375-380.

[3] M. Ahues, Newton methods with Holder derivative, Numerical Functional Analysis and Optimization, 25, N 5-6 (2004), 379-395.

[4] I. K. Argyros, On the convergence of fixed slope iterations, Journal of Mathematics, 38 (2006), 39-44.

[5] I. K. Argyros, Concerning the convergence of a modified Newton-like method, Journal for Analysis and its Applications (ZAA), 18, N 3 (1999), 785-792.

[6] I. K. Argyros, Convergence and Applications of Newton-Type Iterations, Springer, New York, 2008.